\documentclass[11pt]{amsart}

\usepackage{amssymb}
\usepackage[arrow,matrix,curve,cmtip,ps]{xy}
\usepackage{upref}

\newtheorem{thm}{Theorem}[section]   
\newtheorem{cor}[thm]{Corollary}     
\newtheorem{prop}[thm]{Proposition}  

\theoremstyle{definition}

\theoremstyle{remark}
\newtheorem{rem}[thm]{Remark}        

\numberwithin{equation}{section}     

\newcommand{\secref}[1]{Section~\textup{\ref{#1}}}
\newcommand{\thmref}[1]{Theorem~\textup{\ref{#1}}}
\newcommand{\corref}[1]{Corollary~\textup{\ref{#1}}}

\newcommand{\propref}[1]{Proposition~\textup{\ref{#1}}}

\DeclareMathOperator{\ind}{Ind}
\DeclareMathOperator{\id}{id}
\DeclareMathOperator{\ad}{Ad}
\DeclareMathOperator{\aut}{Aut}
\DeclareMathOperator{\Rep}{Rep}

\newcommand{\what}{\widehat}
\newcommand{\wilde}{\widetilde}
\newcommand{\rev}[1]{\widetilde{#1}}
\newcommand{\nonzero}[1]{\quad\text{if $#1$
 \textup(and $0$ if not\textup)}}

\newcommand{\<}{\langle}
\renewcommand{\>}{\rangle}
\newcommand{\dashind}{\text{-}\ind}
\renewcommand{\:}{\colon}
\newcommand{\deltahat}{{\Hat\delta}}
\newcommand{\deltahathat}{\Hat{\Hat\delta}}
\newcommand{\alphahathat}{\Hat{\Hat\alpha}}
\newcommand{\iso}{\xrightarrow{\cong}}
\newcommand{\dec}{^{\operatorname{dec}}}

\renewcommand{\c}[1]{\mathcal #1}

\newcommand{\midtext}[1]{\quad\text{#1}\quad}
\newcommand{\righttext}[1]{\qquad\text{#1 }}

\begin{document}

\title[Three Bimodules]
{Three Bimodules for Mansfield's Imprimitivity Theorem}

\author[Kaliszewski]{S.~Kaliszewski}
\address{Department of Mathematics
\\
Arizona State University
\\
Tempe, AZ 85287}
\email{kaz@math.la.asu.edu}

\author[Quigg]{John Quigg}
\email{quigg@math.la.asu.edu}



\date{January 2, 2000}


\begin{abstract}
For a maximal coaction $\delta$ of a discrete group $G$ on a
$C^*$-algebra $A$ and a normal subgroup $N$ of $G$, 
there are at least three natural $A\times_\delta G\times_{\deltahat|}N
- A\times_{\delta|}G/N$ imprimitivity bimodules: Mansfield's bimodule
$Y_{G/N}^G(A)$; the bimodule assembled by Ng from Green's $A\times_\delta
G\times_{\deltahat} G\times_{\deltahathat|}G/N - A\times_\delta
G\times_{\deltahat|}N$ imprimitivity bimodule $X_N^G(A\times_\delta G)$
and Katayama duality; and the bimodule assembled from
$X_N^G(A\times_\delta G)$ and the crossed-product Mansfield bimodule
$Y_{G/G}^G(A)\times G/N$.  
We show that all three of these are isomorphic, so that the
corresponding inducing maps on representations are identical.  
This can be interpreted as saying that Mansfield and Green
induction are inverses of one another
``modulo Katayama duality''.  These results pass
to twisted coactions; dual results starting with an action are 
also given.
\end{abstract}

\maketitle 


\section{Introduction}
\label{intro-sec}

C.~K.~Ng has recently observed \cite{NgRM} that an abstract Morita
equivalence between a restricted coaction crossed product 
$A\times_{\delta|}
G/N$ and the iterated dual action crossed product $A\times_\delta
G\times_{\hat\delta} N$ can be pieced together from Green's imprimitivity
theorem \cite[Theorem~6]{gre:local}
and Katayama duality \cite[Theorem~8]{kat}, 
thus giving a relatively non-technical, nonconstructive proof of
Mansfield's imprimitivity theorem \cite[Theorem~27]{man}.  
However, in applications --- especially those concerning induced
representations --- it is often necessary to work with an explicit
bimodule.  Because Morita equivalence relations are composed with one
another by tensoring the corresponding imprimitivity bimodules
together, Ng's transitivity argument does implicitly
provide a bimodule.  
Thus the natural question arises as to whether Ng's
bimodule is in fact isomorphic to Mansfield's.  

In more detail: Ng considers a nondegenerate reduced coaction $\delta$ of a
locally compact group $G$ on a $C^*$-algebra $A$, and a closed normal
amenable subgroup $N$ of $G$.  An application of 
Green's theorem to the dual
action $(A\times_\delta G, G, \deltahat)$ gives an $A\times_\delta
G\times_\deltahat G\times_{\deltahathat|}G/N - A\times_\delta
G\times_{\deltahat|}N$ imprimitivity bimodule $X_N^G(A\times_\delta
G)$.  Moreover, looking closely at Katayama's duality theorem, one 
can derive an isomorphism $\Theta$ of $A\times_\delta
G\times_\deltahat G\times_{\deltahathat|}G/N$ onto
$(A\times_{\delta|}G/N)\otimes\c K(L^2(G))$, and this latter algebra is
Morita equivalent to $A\times_{\delta|}G/N$ via the bimodule
$(A\times_{\delta|}G/N)\otimes L^2(G)$. 
Implicitly, then, Ng's $A\times_\delta G\times_{\deltahat|} N -
A\times_{\delta|}G/N$ imprimitivity bimodule is the tensor product
$$\widetilde{X_N^G}(A\times_\delta G)\otimes_\Theta\big(
(A\times_{\delta|}G/N)\otimes L^2(G)\big).$$
(Here the tilde denotes the reverse bimodule.)
Let $Y_{G/N}^G(A)$ denote Mansfield's imprimitivity bimodule.  Then the
question in question is precisely whether
\begin{equation}\label{L2}
X_N^G(A\times_\delta G)\otimes_{A\times G\times N} Y_{G/N}^G(A) \cong
(A\times_{\delta|}G/N)\otimes L^2(G)
\end{equation}
as imprimitivity bimodules.  In other words, modulo crossed-product
duality, are Green and Mansfield induction inverses of one another?

This and related questions concerning actions, twisted actions, and
twisted coactions are addressed in the present paper in the context of
discrete groups and full coactions.  Our approach is to exploit the
natural equivariance of the Mansfield and Green bimodules.  
For instance, in \secref{un}, we consider 
a maximal discrete coaction $(A,G,\delta)$ and Mansfield's 
$A\times_\delta G\times_{\deltahat} G - A$
imprimitivity bimodule $Y_{G/G}^G(A)$,
which carries a $\deltahathat$ -- $\delta$ compatible bimodule
coaction $\delta^Y$.
If in addition $N$ is a normal subgroup of $G$,
\thmref{untw-thm} states that
\begin{equation}\label{Y}
X_N^G(A\times_\delta G)\otimes_{A\times G\times N}
Y_{G/N}^G(A)\cong
Y_{G/G}^G(A)\times_{\delta^Y|}G/N
\end{equation}
as $A\times_\delta G\times_{\deltahat}
G\times_{\deltahathat|} G/N - A\times_{\delta|} G/N$
imprimitivity bimodules.  
The only difference between \eqref{L2} and \eqref{Y} is the bimodule on
the right-hand sides;  but these turn out to be isomorphic
(see the proof of \thmref{Ng-thm}).  Thus Mansfield induction
of representations from $\Rep A\times_{\delta|}G/N$ to $\Rep
A\times_\delta G\times_{\deltahat} N$ via $Y_{G/N}^G(A)$ 
can be ``undone'' by Green
induction via $X_N^G(A\times_\delta G)$ followed by Katayama duality to
get from 
$\Rep A\times_\delta G\times_{\deltahat}G\times_{\deltahathat|}G/N$
back to $\Rep A\times_{\delta|}G/N$. 
In this sense we can usefully view Mansfield and Green induction as
inverse to one another. 
In \secref{dual}, we obtain results dual to those of
\secref{un}, starting with an action  instead of a
coaction.  In \secref{twist} we show that the results of \secref{un}
pass
to twisted coactions, and in \secref{dual twist} we round out this
square of ideas with a set of results for twisted actions.

In \secref{Ng} we return to the comparison between Ng's
bimodule and Mansfield's.  This is done by 
first establishing \eqref{L2} for full coactions and
discrete groups, and then dropping it down to reduced coactions.
In general, we feel that this approach --- 
establishing results first for full coactions, and 
later passing to quotients if results for reduced coactions are
desired --- is more efficient and cleaner conceptually than working
with reduced coactions directly.  
On the other hand, we work with discrete groups simply to avoid many of the
technicalities associated with coactions of general locally compact
groups; also, this is the only context in which we have induced
algebras for coactions, which appear in \secref{dual twist}.  
There is
no reason to believe that the other results in this paper will not hold
in the general case.  In fact, Theorems~\ref{untw-thm} and 
\ref{dual-untw-thm} will
appear in \cite{EKQR} for locally compact groups, but only as a product of
the extensive machinery developed therein.  The proofs here are
considerably more direct, and more instructive.


\section{Preliminaries}
\label{prelim}

In this preliminary section we collect the formulas relevant to 
crossed products and imprimitivity theorems involving actions and 
coactions of \emph{discrete} groups. Because the groups are
discrete, the theory acquires quite an algebraic flavor, and to take 
full advantage of this we translate the standard machinery involving 
locally compact groups to our context.

\subsection*{Coactions}

Let $\delta \: A \to A \otimes C^*(G)$ be a coaction of a discrete
group $G$ on a $C^*$-algebra $A$.  Because $G$ is discrete, the
spectral subspaces $\{A_s : s \in G\}$ of $\delta$ densely span $A$,
and the union $\c A := \bigcup_{s \in G} A_s$
\footnote{more
precisely, the disjoint union, but this abuse will cause no harm,
since the spectral subspaces are linearly independent} forms a Fell
bundle over $G$.  As shown in \cite{eq:induced}, the coaction on $A$
sits ``between'' a ``maximal'' coaction on the full cross-sectional
algebra $C^*(\c A)$ and a ``minimal'' coaction on the reduced
cross-sectional algebra $C^*_r(\c A)$.  Although the crossed product
$A \times_\delta G$ and the covariant representations of the coaction
cannot distinguish among the various possibilities between the
extremes $C^*(\c A)$ and $C^*_r(\c A)$, some other constructions can. 
In particular, the imprimitivity theorems we need require the full
cross-sectional algebra.  Therefore, we \emph{assume throughout that
$A = C^*(\c A)$},
and we call the coaction $\delta$ \emph{maximal} in this case.

The crossed product $A \times_\delta G$ is densely spanned by the
Cartesian product $\c A \times G$, where $A_s \times \{t\}$ has the
obvious vector space structure for all $s,t \in G$, and the
multiplication and involution are given on the generators by
\begin{align*}
(a_r,s)(b_t,u) &= (a_rb_t,u) \nonzero{s=tu}
\\
(a_s,t)^* &= (a_s^*,st).
\end{align*}
Since the coaction $\delta$ is maximal, $A \times_\delta G$ is the
enveloping $C^*$-algebra of the linear span of the generators; that
is, any operation-preserving mapping of the generators into a
$C^*$-algebra $C$ extends uniquely to a homomorphism $A \times_\delta
G \to C$.  If $(B,G,\epsilon)$ is another coaction and $\phi \: A \to
B$ is an equivariant homomorphism (equivalently, $\phi(A_s) \subseteq
B_s$ for each $s \in G$), then we can ``integrate up'' to get a
homomorphism $\phi \times G \: A \times_\delta G \to B \times_\epsilon
G$ defined on the generators by
\[
(\phi \times G)(a_s,t) = (\phi(a_s),t).
\]

The dual action $\hat\delta$ of $G$ on $A \times_\delta G$ is given on
the generators by
\[
\hat\delta_t(a_r,s) = (a_r,st^{-1}).
\]

If $N$ is a normal subgroup of $G$, the coaction $\delta$ restricts to
a maximal coaction $\delta|$ of $G/N$ on $A$, and the crossed product
$A \times_{\delta|} G/N$ is densely spanned by $\c A \times G/N$, with
operations
\begin{align*}
(a_r,sN)(b_t,uN) &= (a_rb_t,uN) \nonzero{sN=tuN}
\\
(a_s,tN)^* &= (a_s^*,stN),
\end{align*}
and maps nondegenerately into $M(A \times_\delta G)$ by
\[
(a_s,tN) \mapsto \sum_{n \in N} (a_s,tn)
\righttext{(strictly convergent).}
\]
There is a \emph{decomposition} coaction $\delta\dec$ of $G$ on the 
restricted crossed product $A \times_{\delta|} G/N$, given on the 
generators by
\[
\delta\dec(a_s,tN) = (a_s,tN) \otimes s.
\]

A coaction $(A,G,\delta)$ is \emph{twisted} over $G/N$ (see 
\cite{pr:twisted}) if there is an orthogonal family $\{p_{tN} : tN \in 
G/N\}$ of projections in $M(A_e)$ which sum strictly to $1$ in
$M(A_e)$, and such that
\[
a_s p_{tN} = p_{stN} a_s
\righttext{for all $s,t\in G$.}
\]
The twisted crossed product $A \times_{\delta,G/N} G$ is the quotient 
of $A \times_\delta G$ by the ideal generated by differences of the 
form
\[
(a_s,t) - (a_s p_{tN},t).
\]
We denote the quotient map by $q \: A \times_\delta G \to A 
\times_{\delta,G/N} G$, and we write
\[
[a_s,t] := q(a_s,t).
\]
The ideal $\ker q$ is invariant under the restriction 
$\hat\delta|_N$, and we
denote the corresponding action of $N$ on $A \times_{\delta,G/N} G$ by 
$\tilde\delta$.
We also define the ``restriction'' $q|$ by the commutative diagram
\[
\xymatrix{
{A \times_{\delta|} G/N}
\ar[r]
\ar[dr]_{q|}
&{M(A \times_\delta G)}
\ar[d]^{q}
\\
&{M(A \times_{\delta,G/N} G).}
}
\]
and we write
\[
[a_s,tN] := q|(a_s,tN).
\]
It is shown in \cite{pr:twisted} that
\[
[a_s,tN] \mapsto a_s p_{tN}
\]
extends to an isomorphism $q(A \times_{\delta|} G/N) \cong A$.

Let $\epsilon$ be a maximal coaction of the quotient $G/N$ on $A$.
It is shown in 
\cite{eq:induced} that there is an \emph{induced} maximal coaction 
$(\ind A,G,\ind\epsilon)$ with spectral subspaces
\[
(\ind A)_s = A_{sN} \times \{s\},
\]
and where the generators have the coordinate-wise operations
\[
(a_{sN},s) (b_{tN},t) = (a_{sN} b_{tN}, st)
\midtext{and} (a_{sN},s)^* = (a_{sN}^*, s^{-1}).
\]

\subsection*{Actions}

Let $\alpha \: G \to \aut B$ be an action of the discrete group $G$ on
a $C^*$-algebra $B$.  The crossed product $B \times_\alpha G$ is
densely spanned by the Cartesian product $B \times G$, where $B \times
\{s\}$ has the obvious vector space structure for all $s \in G$, and
the multiplication and involution are given on the generators by
\[
(a,r)(b,s) = (a\alpha_r(b),rs)
\midtext{and}
(a,r)^* = (\alpha_{r^{-1}}(a^*),r^{-1}).
\]
Again, $B \times_\alpha G$ is the enveloping $C^*$-algebra of the span
of the generators, and if $(C,G,\beta)$ is another action and $\phi \: 
B \to C$ is an equivariant homomorphism, then we can integrate up to 
get a homomorphism $\phi \times G \: B \times_\alpha G \to C 
\times_\beta G$ defined on the generators by
\[
(\phi \times G)(b,s) = (\phi(b),s).
\]

The dual coaction $\hat\alpha$ of $G$ on $B \times_\alpha G$ is given 
on the generators by
\[
\hat\alpha(a,r) = (a,r) \otimes r.
\]

The decomposition action $\alpha\dec$ of $G$ on the restricted crossed 
product $B \times_{\alpha|} N$ is given on the generators by
\[
\alpha\dec_r(a,n) = (\alpha_r(a),rnr^{-1}).
\]

An action $(B,G,\alpha)$ is twisted over $N$ (see \cite{gre:local}) if 
there is a unitary homomorphism $n \mapsto u_n \: N \to M(B)$
such that
\[
\alpha_s(u_n) = u_{sns^{-1}}
\midtext{and}
\alpha_n(b) = u_n b u_n^*.
\]
The twisted crossed product $B \times_{\alpha,N} G$ is the quotient of 
$B \times_\alpha G$ by the ideal generated by differences of the form
\[
(b u_n, s) - (b, ns).
\]
We denote the quotient map by $q \: B \times_\alpha G \to B 
\times_{\alpha,N} G$, and we write
\[
[b,s] := q(b,s).
\]
We denote the dual coaction of $G/N$ on $B \times_{\alpha,N} G$ by 
$\tilde\alpha$.
We also define the ``restriction'' $q|$ by the commutative diagram
\[
\xymatrix{
{B \times_{\alpha|} N}
\ar[r]
\ar[dr]_{q|}
&{B \times_\alpha G}
\ar[d]^{q}
\\
&{B \times_{\alpha,N} G.}
}
\]
and we write
\[
[b,n] := q|(b,n).
\]
It is shown in \cite{gre:local} that
\[
[b,n] \mapsto b u_n
\]
extends to an isomorphism $q(B \times_{\alpha|} N) \cong B$.

Let $\beta$ be an action of the normal subgroup $N$ on $B$.
We identify the induced algebra $\ind B$ as the $c_0$-section 
algebra of a $C^*$-bundle $G \times_N B$ over $G/N$.  Specifically, 
$N$ acts diagonally on the trivial $C^*$-bundle $G \times B \to G$ by 
$n (s,b) := (sn^{-1},nb)$, and the associated orbit space $G \times_N 
B$ has a natural $C^*$-bundle structure over $G/N$: we denote the 
orbit of a pair $(s,b)$ by $[s,b]$, and the fiber over $sN$ is $\{ 
[sn,b] : n \in N, b \in B \}$.  The induced action $(\ind 
B,G,\ind\beta)$ is given on the generators by $\ind\beta_s([t,b]) = 
[st,b]$.

\subsection*{Commutative diagrams of bimodules}

Most of this paper concerns imprimitivity bimodules, but a few times
(in the Sections~\ref{twist} and~\ref{dual twist}) 
we will need the following more general
concept: a \emph{right-Hilbert} $A - B$ bimodule is a Hilbert
$B$-module $X$ equipped with a left $A$-module action by adjointable
maps (which are automatically bounded and $B$-linear).  We will
\emph{assume throughout} that the right inner product is \emph{full}
and the left action is \emph{nondegenerate}.  For example, a
surjective homomorphism $\phi \: A \to B$ determines a right-Hilbert
bimodule ${}_AB_B$ with the obvious right $B$-module action, left
$A$-module action $a \cdot b := \phi(a) b$, and right inner product
$\<b,c\>_B := b^*c$.  Moreover, in this situation any right-Hilbert $B - C$
bimodule $X$ can also be regarded as a right-Hilbert $A - C$ bimodule
with left $A$-module action $a \cdot x := \phi(a) \cdot x$.

We have found it convenient to signify right-Hilbert bimodule
isomorphisms using diagrams: given right-Hilbert bimodules ${}_AX_B$,
${}_BY_C$, and ${}_AZ_C$, when we say the diagram
\[
\xymatrix
{
{A}
\ar[r]^{X}
\ar[dr]_{Z}
&{B}
\ar[d]^{Y}
\\
&{C}
}
\]
commutes, we mean $X \otimes_B Y \cong Z$ as right-Hilbert $A - B$
bimodules, and similarly for rectangular diagrams, {\it etc\/}.  

If ${}_AX_B$ and ${}_CY_D$ are right-Hilbert bimodules and $\phi \: A
\to C$ and $\psi \: B \to D$ are $C^*$-homomorphisms, a linear map
$\Phi \: X \to Y$ is a \emph{right-Hilbert bimodule homomorphism} with
\emph{coefficient maps} $\phi$ and $\psi$ if it preserves the bimodule
actions and the right inner product, that is,
\begin{itemize}
\item[(i)]$\Phi(a \cdot x) = \phi(a) \cdot \Phi(x)$,
\item[(ii)]$\Phi(x \cdot b) = \Phi(x) \cdot \psi(b)$, and
\item[(iii)]$\<\Phi(x),\Phi(y)\>_D = \psi(\<x,y,\>_B)$
\end{itemize}
for all $a\in A$, $b\in B$, and $x,y\in X$.  
If $\phi$ and $\psi$ are isomorphisms then so is
$\Phi$, and if $X$ and $Y$ are imprimitivity bimodules then $\Phi$ is
an imprimitivity bimodule map.

An easy modification of \cite[Lemma 2.5]{sie:morita} shows that the
$B$-linearity condition~(ii) is redundant; in fact, 
if $X_0$ and $A_0$ are dense subspaces of $X$ and $A$, respectively,
and if $\Phi\colon X_0\to Y$ is a linear map satisfying (i) and~(iii) 
above for all $x,y\in X_0$ and $a\in A_0$, then 
$\Phi$ uniquely
extends to a right-Hilbert bimodule homomorphism of ${}_AX_B$ into 
${}_CY_D$.  
Indeed, if $S\subseteq X$ linearly spans $X_0$, $T\subseteq A$ linearly
spans $A_0$, $\Phi$ satisfies (i) and~(iii) on $S$ and $T$ and extends
linearly to $X_0$, then the same conclusion holds.  
We will repeatedly use this fact without comment.

Given a right-Hilbert bimodule homomorphism ${}_AX_B \to {}_CY_D$ with
surjective coefficient maps $\phi \: A \to C$ and $\psi \: B \to D$,
by \cite[Lemma 5.3]{kqr:resind} the diagram
\[
\xymatrix
{
{A}
\ar[r]^{X}
\ar[d]_{\phi}
&{B}
\ar[d]^{\psi}
\\
{C}
\ar[r]_{Y}
&{D}
}
\]
commutes.

\subsection*{Bimodule crossed products}

Let $Z$ be an $A - B$ imprimitivity bimodule, and let $\eta \: Z 
\to Z \otimes C^*(G)$ be a bimodule coaction which is compatible with 
coactions $\delta$ on $A$ and $\epsilon$ on $B$ (see 
\cite{er:multipliers}). Then the spectral subspaces $\{Z_s : s \in 
G\}$ densely span $Z$, and the bimodule crossed product $Z 
\times_\eta G$ is densely spanned by the pairs $\{(x_s,t) : s,t \in 
G, x_s \in Z_s\}$, and is an $A \times_\delta G - B \times_\epsilon 
G$ imprimitivity bimodule with operations given on the generators by
\begin{align*}
(a_r,s)\cdot(x_t,u) &= (a_r\cdot x_t,u) \nonzero{s=tu}
\\
{}_{A\times G}\<(x_r,s),(y_t,u)\> &= ({}_{A}\<x_r,y_t\>,tu) \nonzero{s=u}
\\
(x_r,s)\cdot(b_t,u) &= (x_r\cdot b_t,u) \nonzero{s=tu}
\\
\<(x_r,s),(y_t,u)\>_{B\times G} &= (\<x_r,y_t\>_{B},u)
\nonzero{rs=tu}.
\end{align*}

The dual action $\hat\eta$ of $G$ on $Z \times_\eta G$ is given 
on the generators by
\[
\hat\eta_t(x_r,s) = (x_r,st^{-1}).
\]

Similarly, if $\gamma$ is an action of $G$ on $Z$ which is compatible 
with actions $\alpha$ on $A$ and $\beta$ on $B$, then the crossed 
product $Z \times_\gamma G$ is densely spanned by the Cartesian 
product $Z \times G$, and is an $A \times_\alpha G - B \times_\beta 
G$ imprimitivity bimodule with operations given on the generators by
\begin{align*}
(a,r)\cdot(x,s) &= (a\cdot\alpha_r(x),rs)
\\
{}_{A\times G}\<(x,r),(y,s)\> &= ({}_{A}\<x,\alpha_{rs^{-1}}(y)\>,rs^{-1})
\\
(x,r)\cdot(b,s) &= (x\cdot\alpha_r(b),rs)
\\
\<(x,r),(y,s)\>_{B\times G} &=
(\alpha_{r^{-1}}(\<x,y\>_{B}),r^{-1}s).
\end{align*}

The dual coaction $\hat\gamma$ of $G$ on $Z \times_\gamma G$ is given 
on the generators by
\[
\hat\gamma(x,r) = (x,r) \otimes r.
\]

\subsection*{Imprimitivity theorems}

Let $(A,G,\delta)$ be a maximal discrete coaction, and let $N$ be a
normal subgroup of $G$.  By the version of Mansfield's imprimitivity
theorem due to Echterhoff and the second author
\cite[Theorem~3.1]{eq:full}, there exists an $A \times_\delta G
\times_{\hat\delta|} N - A \times_{\delta|} G/N$ imprimitivity
bimodule $Y_{G/N}^G(A)$.  Mansfield's bimodule is densely spanned by
the Cartesian product $\c A \times G$, with operations given on the
generators by
\begin{align*}
(a_r,s,n) \cdot (b_t,u) &= (a_rb_t,un^{-1}) \nonzero{sn=tu}
\\
{}_{A \times G \times N}\<(a_r,s),(b_t,u)\> &= (a_rb_t^*,ts,s^{-1}u)
\nonzero{sN=uN}
\\
(a_r,s) \cdot (b_t,uN) &= (a_rb_t,t^{-1}s) \nonzero{sN=tuN}
\\
\<(a_r,s),(b_t,u)\>_{A \times G/N} &= (a_r^*b_t,uN) \nonzero{rs=tu}
\end{align*}
It is easy to see that $Y(A)$ is functorial in the sense that if
$(B,G,\epsilon)$ is another coaction and $\phi \: A \to B$ is an
equivariant homomorphism then $(a_s,t) \mapsto (\phi(a_{s}),t)$
extends to an imprimitivity bimodule homomorphism $Y(A) \to Y(B)$ with
coefficient homomorphisms $\phi \times G \times N$ and $\phi \times
G/N$.

Dually, for an action $(B,G,\alpha)$ and a normal subgroup $N$ of $G$,
Green's imprimitivity theorem \cite[Proposition~3]{gre:local} provides
a $B \times_\alpha G \times_{\hat\alpha|} G/N - B \times_{\alpha|} N$
imprimitivity bimodule $X_N^G(B)$.  Green's bimodule is densely
spanned by the Cartesian product $B \times G$, with operations given
on the generators by
\begin{align*}
(a,r,sN) \cdot (b,t) &= (a\alpha_r(b),rt) \nonzero{sN=tN}
\\
{}_{B \times G \times G/N}\<(a,r),(b,s)\> &=
(a \alpha_{rs^{-1}}(b^*),rs^{-1},sN)
\\
(a,r) \cdot (b,n) &= (a \alpha_r(b),rn)
\\
\<(a,r),(b,s)\>_{B \times N} &= (\alpha_{r^{-1}}(a^*b),r^{-1}s)
\nonzero{rN=sN}
\end{align*}
$X(B)$ is functorial in the sense that if $(C,G,\beta)$ is another
action and $\phi \: B \to C$ is an equivariant homomorphism then
$(b,s) \mapsto (\phi(b),s)$ extends to an imprimitivity bimodule
homomorphism $X(B) \to X(C)$ with coefficient homomorphisms $\phi
\times G \times G/N$ and $\phi \times N$.

If the coaction $(A,G,\delta)$ is twisted over $G/N$, then the
quotient map $q \: A \times_{\delta} G \to A \times_{\delta,G/N} G$
restricts to a surjection $q| \: A \times_{\delta|} G/N \to A$, giving
us an ideal $\ker q|$ of $A \times G/N$.  Inducing across Mansfield's
bimodule $Y$ via the Rieffel correspondence gives an ideal $Y \dashind
(\ker q|)$ of $A \times_\delta G \times_{\hat\delta|} N$, and
\cite[Theorem 4.1]{pr:twisted} shows that this ideal is precisely
$\ker (q \times N)$.  Rieffel's theory thus gives an $A
\times_{\delta,G/N} G \times_{\tilde\delta} N - A$ imprimitivity
bimodule
\[
Z_{G/N}^G(A) := Y/(Y \cdot \ker q|).
\]
Moreover, the diagram
\[
\xymatrix
{
{A \times_\delta G \times_{\hat\delta|} N}
\ar[r]^-{Y}
\ar[d]_{q \times N}
&{A \times_{\delta|} G/N}
\ar[d]^{q|}
\\
{A \times_{\delta,G/N} G \times_{\tilde\delta} N}
\ar[r]_-{Z}
&{A}
}
\]
commutes.

Dually, if the action $(B,G,\alpha)$ is twisted over $N$, by 
\cite[Corollary 5]{gre:local} we have a $B 
\times_{\alpha,N} G \times_{\hat\alpha|} G/N - B$ imprimitivity 
bimodule
\[
W_N^G(B) := X/(X \cdot \ker q|),
\]
and we get a commutative diagram
\[
\xymatrix
{
{B \times_\alpha G \times_{\hat\alpha|} G/N}
\ar[r]^-{X}
\ar[d]_{q \times G/N}
&{B \times_{\alpha|} N}
\ar[d]^{q|}
\\
{B \times_{\alpha,N} G \times_{\tilde\alpha} G/N}
\ar[r]_-{W}
&{B.}
}
\]

If $\epsilon$ is a maximal coaction of the quotient $G/N$ on $A$, the
recent imprimitivity theorem for induced coactions \cite[Theorem
4.1]{eq:induced} gives an $\ind A \times_{\ind\epsilon} G - A
\times_\epsilon G/N$ imprimitivity bimodule $U = U(A)$ densely spanned
by the subset
\[
\{ (a_{sN},t) : sN \in G/N, a_{sN} \in A_{sN}, t \in G \}
\]
of the Cartesian product $A \times G$, with operations given on the 
generators by
\begin{align*}
(a_{sN},t,r) (b_{uN},v) &= (a_{sN} b_{uN}, tv) \nonzero{r = v}
\\
(a_{sN},t) (b_{uN},vN) &= (a_{sN} b_{uN}, t) \nonzero{tN = suvN}
\\
{}_{\ind A \times G} \< (a_{sN},t), (b_{uN},v) \>
&= (a_{sN} b_{uN}^*, tv^{-1}, v) \nonzero{s^{-1}tN = u^{-1}vN}
\\
\< (a_{sN},t), (b_{uN},v) \>_{A \times G/N}
&= (a_{sN}^* b_{uN}, u^{-1}vN) \nonzero{t = v}.
\end{align*}
$U(A)$ is functorial in the sense that if $(B,G/N,\eta)$ is another
coaction and $\phi \: A \to B$ is an equivariant homomorphism then
$(a_sN,t) \mapsto (\phi(a_{sN}),t)$ extends to an imprimitivity
bimodule homomorphism $U(A) \to U(B)$ with coefficient homomorphisms
$\ind \phi \times G$ and $\phi \times G/N$.

Dually, if $\beta$ is an action of the normal subgroup $N$ on $B$, the
imprimitivity theorem for induced actions (sometimes attributed to
Green \cite[Theorem 17]{gre:local}) gives an $\ind B
\times_{\ind\beta} G - B \times_\beta N$ imprimitivity bimodule $V =
V(B)$ densely spanned by the Cartesian product $B \times G$, with
operations given on the generators by
\begin{align*}
([tr,b],t) (c,r) &= (bc,tr)
\\
(b,s) (c,n) &= (\beta_{n^{-1}}(bc),sn)
\\
{}_{\ind B \times G} \< (b,s), (c,t) \> &= ([s,bc^*],st^{-1})
\\
\< (b,s), (c,sh) \>_{B \times N} &= (b^* \beta_h(c),h)
\end{align*}
$V(B)$ is functorial in the sense that if $(C,N,\beta)$ is another
action and $\phi \: B \to C$ is an equivariant homomorphism then
$(b,s) \mapsto (\phi(b),s)$ extends to an imprimitivity bimodule
homomorphism $V(B) \to V(C)$ with coefficient homomorphisms $\ind \phi
\times G$ and $\phi \times N$.


\section{The Mansfield-Green Triangle}
\label{un}

In this section we show a curious duality between the Mansfield and 
Green imprimitivity theorems.  \thmref{untw-thm} will show 
that, roughly speaking, and modulo crossed product duality, Mansfield
and Green induction are inverse processes.

Let $(A,G,\delta)$ be a maximal discrete coaction, and let $N$ be a 
normal subgroup of $G$.  Not only do we have Mansfield's 
$A\times_\delta G \times_{\hat\delta|} N -
A\times_{\delta|}G/N$ imprimitivity bimodule $Y_{G/N}^G(A)$, but 
also, replacing $N$ by $G$, an $A \times_\delta G \times_{\hat\delta} 
G - A$ imprimitivity bimodule $Y_{G/G}^G(A)$.  There is a 
$\deltahathat - \delta$ compatible coaction $\delta_Y$ of $G$ on 
$Y_{G/G}^G(A)$ (\cite[Remark~3.2]{eq:full}) determined by
\[
\delta_Y(a_r,s) = (a_r,s)\otimes s^{-1}.  
\]

\begin{thm}
\label{untw-thm}
Let $(A,G,\delta)$ be a maximal coaction and let $N$ be a normal
subgroup of $G$.  Then the diagram
\[
\xymatrix
@C+2pt
@L+3pt
{
{A \times_\delta G \times_{\hat\delta|} N}
\ar[r]^-{Y_{G/N}^G(A)}
&{A \times_{\delta|} G/N} 
\\
{A \times_\delta G \times_{\hat\delta} G \times_{\deltahathat\bigm|} G/N.}
\ar[u]^{X_N^G(A\times G)}
\ar[ur]_<(.6){Y_{G/G}^G(A)\times G/N}
}
\]
commutes.
\end{thm}

\begin{proof}
There is an imprimitivity bimodule isomorphism
\[
\Phi \: Y(A) \times G/N \otimes \wilde{Y(A)}
\iso X(A \times G)
\]
defined on the generators by
\[
\Phi((a_r,s,tN) \otimes \rev{(b_u,v)}) =
(a_r b_u^*,us,s^{-1}v) \nonzero{tN=vN},
\]
since straightforward calculations verify that the above mapping on
generators preserves the left action and the right inner
product.
\end{proof}

\begin{rem}
\label{untw-rem}
To motivate the formula for $\Phi$, note that 
$$\Phi((a_r,s,tN)\otimes\rev{(c_u,v)})
= {}_{A\times G\times G}\<(a_r,s),(c_u,v)\>\nonzero{tN=vN},$$
where $(a_r,s)$ and $(c_u,v)$ are viewed as elements of $Y_{G/G}^G(A)$,
and the inner product is viewed as taking values in
$X_N^G(A\times_\delta G)$.  
\end{rem}

The next two results will not be needed until \secref{twist}; we
include them here because they don't involve twists, and are of general
interest. 

\begin{prop}
\label{move-G/N-prop}
Let $(A,G,\delta)$ be a coaction, and let $N$ be a normal subgroup of
$G$.  Then
\[
A\times_\delta
G\times_{\hat\delta}G\times_{\deltahathat\bigm|}G/N \cong
A\times_{\delta|}G/N \times_{\delta\dec}G
\times_{\what{\delta\dec}}G.
\]
\end{prop}

\begin{proof}
Straightforward calculations verify that the mapping
\[
(a_r,s,t,uN) \mapsto (a_r,stuN,s,t)
\]
on the generators preserves the operations, hence extends to a 
$C^*$-isomorphism.
\end{proof}

\begin{prop}
\label{G/N-bimod-prop}
Let $(A,G,\delta)$ be a coaction and $N$ a normal subgroup of
$G$. 
Then the diagram
\[
\xymatrix{
&{A \times_{\delta|} G/N}
\\
{A \times_\delta G \times_{\hat\delta} G \times_{\deltahathat\bigm|} 
G/N}
\ar[ur]^{Y_{G/G}^G(A)\times G/N}
\ar[r]^-{\cong}
&{A\times_{\delta|}G/N \times_{\delta\dec}G
\times_{\what{\delta\dec}}G.}
\ar[u]_{Y_{G/G}^G(A\times G/N)}
}
\]
commutes, where the isomorphism is that of \propref{move-G/N-prop}.
\end{prop}

\begin{proof}
There is an imprimitivity bimodule isomorphism
\[
\Theta \: Y(A) \times G/N \iso Y(A \times G/N)
\]
defined on the generators by
\[
\Theta(a_r,s,tN) = (a_r,tN,s),
\]
since straightforward calculations verify that the above mapping on
generators preserves the left action and the right inner product.
\end{proof}


\section{The Dual Triangle}
\label{dual}

The results of this section are dual to those of the previous
section,  in the sense that actions correspond to coactions, Green
bimodules correspond to Mansfield bimodules, and subgroups correspond
to quotient groups.
The only additional apparatus we need is to
observe that if $(B,G,\alpha)$ is a discrete action then
there is an $\alphahathat - \alpha$ compatible
action $\alpha^X$ of $G$ on $X_e^G(B)$ given on the generators by 
\[
\alpha_r^X(a,s) = (a,sr^{-1}).
\]

\begin{thm}
\label{dual-untw-thm}
Let $(B,G,\alpha)$ be an action an let $N$ be a normal subgroup
of $G$.  Then the diagram
\[
\xymatrix
@C+20pt
@L+3pt
{
{B\times_\alpha G \times_{\hat\alpha} G \times_{\alphahathat|} N}
\ar[r]^-{Y_{G/N}^G(B\times G)}
\ar[dr]_{X_e^G(B)\times N}
&{B\times_\alpha G \times_{\hat\alpha|} G/N}
\ar[d]^{X_N^G(B)}
\\
&{B\times_{\alpha|} N.}
}
\]
commutes.
\end{thm}

\begin{proof}
There is an imprimitivity bimodule isomorphism
\[
\Phi \: (X_e^G(B)\times N)\otimes \widetilde{X_N^G}(B)
\iso Y_{G/N}^G(B\times G)
\]
defined on the generators by
\[
\Phi((a,r,n)\otimes\rev{(b,s)}) =
(a\alpha_{rns^{-1}}(b^*),rns^{-1},sn^{-1}),
\]
since straightforward calculations verify that the above mapping on
generators preserves the left action and the right inner product.
\end{proof}

\begin{rem}
\label{dual-untw-rem}
To motivate the formula for $\Phi$, note that 
$$\Phi((a,r,n)\otimes\rev{(b,s)}) =
{}_{B\times G\times G}\<(a,r),\alpha_n^X(b,s)\>,$$
where $(a,r)$ and $(b,s)$ are viewed as elements of $X_e^G(B)$, and the
inner product is viewed as taking values in $Y_{G/N}^G(B\times_\alpha
G)$.  
\end{rem}

In analogy with the previous section, the next two results will not be
needed until \secref{dual twist}; they are presented here for
convenience and general interest. 

\begin{prop}
\label{move-N-prop}
Let $(B,G,\alpha)$ be an action and let $N$ be a normal subgroup of
$G$.  Then 
\[
B\times_\alpha G\times_{\hat\alpha} G\times_{\alphahathat|}N \cong
B\times_{\alpha|}N\times_{\alpha\dec}G\times_{\what{\alpha\dec}}G.
\]
\end{prop}

\begin{proof}
Straightforward calculations verify that the mapping
\[
(a,r,s,n) \mapsto (a,rsns^{-1}r^{-1},rsn^{-1}s^{-1},sn)
\]
on the generators preserves the operations, hence extends to a 
$C^*$-isomorphism.
\end{proof}

\begin{prop}
\label{N-bimod-prop}
Let $(B,G,\alpha)$ be an action and let $N$ be a normal subgroup of
$G$.  Then the diagram
\[
\xymatrix
@L+2pt
@C+15pt
{
{B\times_\alpha G \times_{\hat\alpha} G \times_{\alphahathat|} N}
\ar[d]_{\cong}
\ar[dr]^<(.5){X_e^G(B)\times N}
\\
{B\times_{\alpha|} N \times_{\alpha\dec} G \times_{\what{\alpha\dec}}
G}
\ar[r]_-{X_e^G(B\times N)}
&{B\times_{\alpha|} N.}
}
\]
commutes.
\end{prop}

\begin{proof}
There is an imprimitivity bimodule isomorphism
\[
\Theta \:
X(B) \times N
\iso
X(B \times N)
\]
defined on the generators by
\[
\Theta(a,r,n) = (a,rnr^{-1},r),
\]
since straightforward calculations verify that the above mapping on
generators preserves the left action and the right inner product.
\end{proof}


\section{The twisted Mansfield-Green square}
\label{twist}

Let $(A,G,\delta)$ be a maximal discrete coaction, and let $N$ be a 
normal subgroup of $G$.
Combining \thmref{untw-thm} and \corref{G/N-bimod-prop}, we get a 
commutative rectangle
\begin{equation}
\label{untw sq}
\xymatrix
@C+2pt
@L+3pt
{
{A \times_\delta G \times_{\hat\delta|} N}
\ar[r]^-{Y_{G/N}^G(A)}
&{A \times_{\delta|} G/N} 
\\
{A \times_\delta G \times_{\hat\delta} G \times_{\deltahathat\bigm|} G/N}
\ar[u]^{X_N^G(A \times G)}
\ar[r]^-{\cong}_-{\Upsilon}
&{A \times_{\delta|} G/N \times_{\delta\dec} G
\times_{\what{\delta\dec}} G.}
\ar[u]_{Y_{G/G}^G(A \times G/N)}
}
\end{equation}
Now suppose the coaction $\delta$ is twisted over $G/N$.  Then the top
arrow of the diagram \eqref{untw sq} has
\[
\xymatrix
{
{A \times_{\delta,G/N} G \times_{\tilde\delta} N}
\ar[r]^-{W}
&{A}
}
\]
as a quotient.  The imprimitivity bimodules in \eqref{untw sq} above
determine corresponding ideals of the bottom corners, and we can form
a quotient commutative rectangle with upper right corner $A$.  What
happens to the rest of the diagram?  We will answer this question in
the present section.

However, we first modify the lower left corner of the diagram
\eqref{untw sq}; the action $\tilde\delta$ of $N$ on $A \times_{G/N}
G$ does not extend to $G$, so the Green bimodule $X_N^G$ on the left
edge of \eqref{untw sq} will not pass to a Green bimodule in the
quotient.  Rather, it will be more appropriate to use the bimodule
arising from the imprimitivity theorem for induced actions.
The action $\tilde\delta$ of $N$ on 
$A \times_{G/N} G$ induces to an action $\ind \tilde\delta$ of
$G$ on the induced algebra $\ind (A \times_{G/N} G)$, and we
have an $\ind (A \times_{G/N} G) \times_{\ind \tilde\delta} G -
(A \times_{G/N} G) \times_{\tilde\delta} N$ imprimitivity
bimodule $V_N^G(A \times_{G/N} G)$.

It is shown in \cite[Theorem 4.4]{qr:induced} that
\[
(a_s p_{tN},r) \mapsto \bigl[ r^{-1}t,[a_s,t] \bigr]
\]
extends to an isomorphism $A \times_\delta G \cong \ind (A
\times_{\delta,G/N} G)$ which is equivariant for the actions
$\hat\delta$ and $\ind\tilde\delta$ of $G$.  Turning this around and
integrating up, we get an isomorphism
\[
\ind (A \times_{\delta,G/N} G) \times_{\ind\tilde\delta} G
\cong
A \times_\delta G \times_{\hat\delta} G.
\]

\begin{thm}  
If $(A,G,\delta)$ is a maximal discrete coaction which is twisted 
over $G/N$, the diagram
\[
\xymatrix{
{A \times_{\delta,G/N} G \times_{\tilde\delta} N}
\ar[r]^-{Z_{G/N}^G(A)}
&{A}
\\
{\ind (A \times_{\delta,G/N} G) \times_{\ind\tilde\delta} G}
\ar[u]^{V_N^G(A \times_{G/N} G)}
\ar[r]_-{\cong}
&{A \times_\delta G \times_{\hat\delta} G}
\ar[u]_{Y_{G/G}^G(A)}
}
\]
commutes.
\end{thm}

\begin{proof}
The desired diagram is the inner rectangle of the diagram
\begin{equation}
\label{untw-tw}
\xymatrix
{
{A \times G \times N}
\ar[rrr]^{Y_{G/N}^G(A)}
\ar[dr]|{q \times N}
&&&{A \times G/N}
\ar[dl]|{q|}
\\
&{A \times_{G/N} G \times N}
\ar[r]^-{Z_{G/N}^G(A)}
&{A}
\\
&{\ind (A \times_{G/N} G) \times G}
\ar[u]|{V_N^G(A \times_{G/N} G)}
\ar[r]_-{\cong}
&{A \times G \times G}
\ar[u]|{Y_{G/G}^G(A)}
\\
{\ind (A \times G) \times G}
\ar[uuu]|{V_N^G(A \times G)}
\ar[ur]|{\ind q \times G}
\ar@{-->}[rrr]_{\cong}
&&&{A \times G/N \times G \times G.}
\ar[ul]|{q| \times G \times G}
\ar[uuu]|{Y_{G/G}^G(A \times G/N)}
}
\end{equation}
(Here and in Diagram~\eqref{V-X} the action and coaction symbols 
have been omitted for clarity.)
We will show how to fill in the bottom arrow so that each of the
outer rectangle and the 
top, bottom, left, and right quadrilaterals commute.
Since $\ind q \times G$ is surjective, the result will then follow 
from standard bimodule techniques.

Consider the diagram
\begin{equation}
\label{V-X}
\xymatrix
{
{A \times G \times N}
\ar[rrr]^-{Y_{G/N}^G(A)}
&&&{A \times G/N}
\\
\\
&{A \times G \times G \times G/N}
\ar[uul]_{X_N^G(A\times G)}
\ar[rrd]^{\Upsilon}_{\cong}
\\
{\ind (A \times G) \times G}
\ar[uuu]^{V_N^G(A\times G)}
\ar@{-->}[ur]^{\cong}
\ar@{-->}[rrr]_{\cong}
&&&{A \times G/N \times G \times G.}
\ar[uuu]_{Y_{G/G}^G(A\times G/N)}
}
\end{equation}
Since the action $\hat\delta|$ of $N$ extends to $G$, it follows from 
standard facts concerning induced actions that
\[
([s,a_t,r],u) \mapsto
(a_t,rs^{-1},u,u^{-1}sN)
\]
extends to an isomorphism 
$\ind (A \times_\delta G) \times_{\ind\hat\delta|} G
\cong A \times_\delta G \times_{\hat\delta} G \times_{\deltahathat\bigm|} G/N$. 
Then an easy check on the 
generators shows that
\[
(a_s,t,r) \mapsto (a_s,tr^{-1},r)
\]
extends to an isomorphism $V(A \times_\delta G) \cong X(A \times_\delta G)$ of 
$\ind (A \times_\delta G) \times_{\ind\hat\delta|} G 
- A \times_\delta G \times_{\hat\delta|} N$ imprimitivity 
bimodules. This shows the left triangle of the diagram \eqref{V-X} 
commutes. The inner quadrilateral in \eqref{V-X} is 
the commutative diagram \eqref{untw sq}. We \emph{define}
the isomorphism $\ind (A \times_\delta G) \times_{\ind\hat\delta|} G 
\cong A \times_{\delta|} G/N \times_{\delta\dec} 
G \times_{\what{\delta\dec}} G$ at the bottom arrow of 
\eqref{V-X} so that the bottom 
triangle commutes. On the generators, this isomorphism is given by
\[
([s,a_t,r],u) \mapsto (a_t,rN,rs^{-1},u).
\]
Thus the outer rectangle in the diagram \eqref{untw-tw} commutes.

We noticed in \secref{prelim} that the
top quadrilateral in \eqref{untw-tw} commutes.

For the right quadrilateral, $Y(B)$ is functorial in $B$,
so the homomorphism $q| \: A \times_{\delta|} G/N 
\to A$ yields an imprimitivity bimodule homomorphism $Y(q|) \: Y(A 
\times_{\delta|} G/N) \to Y(A)$ with the desired coefficient homomorphisms.  
By 
\cite[Lemma 5.3]{kqr:resind}, this implies the quadrilateral commutes.

Similarly, the left quadrilateral commutes by functoriality 
of $V(B)$: the homomorphism $q \: A \times_\delta G \to 
A \times_{\delta,G/N} G$ 
yields an imprimitivity bimodule homomorphism $V(q) \: V(A \times_\delta G) 
\to V(A \times_{\delta,G/N} G)$ with the desired coefficient homomorphisms.

Finally, the bottom quadrilateral in \eqref{untw-tw} commutes by a 
routine computation on the generators.
\end{proof}


\section{The twisted dual square}
\label{dual twist}

In this section we introduce a twist into \thmref{dual-untw-thm}, just 
as in the preceding section we threw a twist into \thmref{untw-thm}; 
unsurprisingly, the development will closely parallel that of 
\secref{twist}.

Let $(B,G,\alpha)$ be a discrete action which is twisted over a normal 
subgroup $N$ in the sense of \cite{gre:local}.
\thmref{dual-untw-thm} and \corref{N-bimod-prop} together give a 
commutative rectangle
\begin{equation}
\label{dual untw sq}
\xymatrix
@C+2pt
@L+3pt
{
{B \times_\alpha G \times_{\hat\alpha|} G/N}
\ar[r]^-{X_{N}^G(B)}
&{B \times_{\alpha|} N} 
\\
{B \times_\alpha G \times_{\hat\alpha} G \times_{\alphahathat|} N}
\ar[u]^{Y_{G/N}^G(B \times G)}
\ar[r]^-{\cong}_-{\Upsilon}
&{B \times_{\alpha|}N \times_{\alpha\dec}G
\times_{\what{\alpha\dec}}G.}
\ar[u]_{X_{\{e\}}^G(B \times N)}
}
\end{equation}
As in the preceding section, in order to form a suitable quotient 
diagram we need to replace the lower left corner by an induced algebra.

The dual coaction $\tilde\alpha$ of $G/N$ on the twisted crossed 
product $B \times_{\alpha,N} G$ induces to a coaction $\ind 
\tilde\alpha$ of $G$ on the induced algebra $\ind (B 
\times_{\alpha,N} G)$, and we have an $\ind (B \times_{\alpha,N} G) 
\times_{\ind \tilde\alpha} G - (B \times_{\alpha,N} G) 
\times_{\tilde\alpha} G/N$ imprimitivity bimodule $U_{G/N}^G(B 
\times_{\alpha,N} G)$.

It is shown in \cite[Theorem 5.6]{eq:induced} that
\[
(b,s) \mapsto ([b,s],s)
\]
extends to an isomorphism
$B \times_\alpha G \cong \ind (B \times_{\alpha,N} G)$
which is equivariant for the coactions $\hat\alpha$ and 
$\ind\wilde\alpha$ of $G$.
Turning this around and integrating up, we 
get an isomorphism
\[
\ind (B \times_{\alpha,N} G) \times_{\ind\tilde\alpha} G
\cong
B \times_\alpha G \times_{\hat\alpha} G.
\]

\begin{thm}

If $(B,G,\alpha)$ is a discrete action which is twisted 
over $N$, the diagram
\[
\xymatrix{
{B \times_{\alpha,N} G \times_{\tilde\alpha} G/N}
\ar[r]^-{W_{N}^G(B)}
&{B}
\\
{\ind (B \times_{\alpha,N} G) \times_{\ind\tilde\alpha} G}
\ar[u]^{U_{G/N}^G(B \times_{\alpha,N} G)}
\ar[r]_-{\cong}
&{B \times_\alpha G \times_{\hat\alpha} G}
\ar[u]_{X_{\{e\}}^G(B)}
}
\]
commutes.
\end{thm}

\begin{proof}
The desired diagram is the inner rectangle of the diagram
\begin{equation}
\label{dual untw-tw}
\xymatrix
{
{B \times G \times G/N}
\ar[rrr]^{X_{N}^G(B)}
\ar[dr]|{q \times G/N}
&&&{B \times N}
\ar[dl]|{q|}
\\
&{B \times_{N} G \times G/N}
\ar[r]^-{W_{N}^G(B)}
&{B}
\\
&{\ind (B \times_{N} G) \times G}
\ar[u]|{U_N^G(B \times_{N} G)}
\ar[r]_-{\cong}
&{B \times G \times G}
\ar[u]|{X_{\{e\}}^G(B)}
\\
{\ind (B \times G) \times G}
\ar[uuu]|{U_N^G(B \times G)}
\ar[ur]|{\ind q \times G}
\ar@{-->}[rrr]_{\cong}
&&&{B \times N \times G \times G.}
\ar[ul]|{q| \times G \times G}
\ar[uuu]|{X_{\{e\}}^G(B \times N)}
}
\end{equation}

Consider the diagram
\begin{equation}
\label{U-Y}
\xymatrix
{
{B \times G \times G/N}
\ar[rrr]^-{X_N^G(B)}
&&&{B \times N}
\\
\\
&{B \times G \times G \times N}
\ar[uul]_{Y_{G/N}^G(B\times G)}
\ar[rrd]^{\Upsilon}_{\cong}
\\
{\ind (B \times G) \times G}
\ar[uuu]^{U_{G/N}^G(B\times G)}
\ar@{-->}[ur]^{\cong}
\ar@{-->}[rrr]_{\cong}
&&&{B \times N \times G \times G.}
\ar[uuu]_{X_{\{e\}}^G(B\times N)}
}
\end{equation}
It follows from 
\cite[Remark 3.3]{eq:full} that the map
\[
(b,s,sn,t) \mapsto (b,s,nt,t^{-1}n^{-1}t)
\]
extends to an isomorphism 
$\ind (B \times_\alpha G) \times_{\ind\hat\alpha|} G \cong 
B \times_\alpha G \times_{\hat\alpha} G \times_{\alphahathat|} N$, 
and this serves as the 
left-hand coefficient map for an isomorphism 
$U(B \times_\alpha G) \cong Y(B \times_\alpha G)$, hence
the left triangle of the diagram \eqref{U-Y} 
commutes. The inner quadrilateral is 
the commutative diagram \eqref{dual untw sq}. We define
the isomorphism 
$\ind (B \times_\alpha G) \times_{\ind\hat\alpha|} G \cong 
B \times_{\alpha|} N \times_{\alpha\dec} G \times_{\what{\alpha\dec}} G$ 
at the bottom arrow of \eqref{U-Y} so that the bottom 
triangle commutes. On the generators, this isomorphism is given by
\[
(b,s,t,r) \mapsto (b,st^{-1},t,r).
\]
Thus the outer rectangle in the diagram \eqref{dual untw-tw} commutes.

We noticed in \secref{prelim} that the
top quadrilateral in \eqref{dual untw-tw} commutes.

$X(A)$ is functorial in $A$,
so the homomorphism $q| \: B \times_{\alpha|} N 
\to B$ yields an imprimitivity bimodule homomorphism $X(q|) \: X(B 
\times_{\alpha|} N) \to X(B)$ with the desired coefficient homomorphisms, so
the right quadrilateral commutes.

Similarly, the left quadrilateral commutes because by functoriality 
of $U(A)$ the homomorphism 
$q \: B \times_\alpha G \to B \times_{\alpha,N} G$ 
yields an imprimitivity bimodule homomorphism 
$U(q) \: U(B \times_\alpha G) \to U(B \times_{\alpha,N} G)$ 
with the desired coefficient homomorphisms.

Finally, the bottom quadrilateral in \eqref{dual untw-tw} commutes by a 
routine computation on the generators.
\end{proof}


\section{Ng's Bimodule}
\label{Ng}

We now return to the comparison between Ng's bimodule and Mansfield's,
beginning with maximal coactions and full crossed products.  In this
context, by ``Ng's bimodule'' we mean the bimodule gotten from the
lower three sides of Diagram~\eqref{ng diag}; the map $\Theta$ will be
defined in the proof of \thmref{Ng-thm} by a construction parallel to Ng's.

\begin{thm}
\label{Ng-thm}
If $(A,G,\delta)$ is a maximal coaction of a discrete group $G$ and
$N$ is a normal subgroup of $G$, then Ng's bimodule is isomorphic to
Mansfield's; that is, the diagram
\begin{equation}
\label{ng diag}
\xymatrix{
{A \times_\delta G \times_{\hat\delta|} N}
\ar[r]^-{Y_{G/N}^G(A)}
&{A \times_{\delta|} G/N}
\\
{A \times_\delta G \times_{\hat\delta} G \times_{\deltahathat\bigm|} G/N}
\ar[u]^{X_N^G(A\times G)}
\ar[r]^-{\cong}_-{\Theta}
&{(A \times_{\delta|} G/N) \otimes \c K(\ell^2(G))}
\ar[u]_{(A \times G/N) \otimes \ell^2(G)}
}
\end{equation}
commutes.
\end{thm}

\begin{proof}
The desired diagram is the outer rectangle of
\[
\xymatrix
@L+2pt
{
{A \times_\delta G \times_{\hat\delta|} N}
\ar[r]^-{Y_{G/N}^G(A)}
&{A \times_{\delta|} G/N}
\\
{A \times_\delta G \times_{\hat\delta} G \times_{\deltahathat\bigm|} G/N}
\ar[u]^{X_N^G(A\times G)}
\ar[r]_-{\Theta}
\ar[ur]|(.5){Y_{G/G}^G(A) \times G/N}
&{(A \times_{\delta|} G/N) \otimes \c K(\ell^2(G)).}
\ar[u]_{(A \times G/N) \otimes \ell^2(G)}
}
\]
The upper left triangle commutes by \thmref{untw-thm}, so we must show 
the lower right triangle commutes.

We construct the isomorphism $\Theta$ as a composition
\begin{multline}\label{theta-eq}
A \times_\delta G \times_{\hat\delta} G \times_{\deltahathat\bigm|} G/N
\xrightarrow{\Theta_1}
(A \otimes \c K) \times_{\epsilon_1} G/N
\xrightarrow{\Theta_2}
(A \otimes \c K) \times_{\epsilon_2} G/N
\\
\xrightarrow{\Theta_3}
(A \times_{\delta|} G/N) \otimes \c K.
\end{multline}
Here, $\epsilon_2$ is the coaction
\[
\epsilon_2 := (\id \otimes \sigma) \circ (\delta| \otimes \id)
\]
of $G/N$ on $A \otimes \c K$,
where $\sigma \: A \otimes C^*(G/N) \otimes \c K \to A \otimes \c K
\otimes C^*(G/N)$ is the flip isomorphism.
It follows from \cite[Lemma 1.16 (b)]{qui:fullreduced} (see also 
\cite{rae:representation}) that there is an
isomorphism $\Theta_3$ of $(A \otimes \c K) \times_{\epsilon_2} G/N$
onto $(A \times_{\delta|} G/N) \otimes \c K$ defined on the generators by
\[
\Theta_3(a_{sN} \otimes b,tN) = (a_{sN},tN) \otimes b.
\]

Perpetuating our perverse numbering scheme, we use
$\epsilon_2$ to define the coaction
\[
\epsilon_1 := \ad (1 \otimes U) \circ \epsilon_2
\]
of $G/N$ on $A \otimes \c K$, where $U$ is the unitary element of 
$M(c_0(G) \otimes C^*(G/N)) \subseteq M(\c K \otimes C^*(G/N))$
determined by the bounded
function $U(s) = sN$ from $G$ to $C^*(G/N)$. It is easy to see that 
$U$ is an $\epsilon_1$-cocycle (more precisely, the obvious analogue 
for full coactions of the more usual cocycles for reduced coactions---see 
\cite{lprs}). It follows from
\cite[Theorem 2.9]{lprs}
(also see \cite[Proposition 2.8]{qr:induced})
that there is an isomorphism $\Theta_2$ of $(A \otimes \c K)
\times_{\epsilon_1} G/N$ onto $(A \otimes \c K)
\times_{\epsilon_2} G/N$ defined by
\[
\Pi_2 \circ \Theta_2
= \ad (\id_{A \otimes \c K} \otimes \lambda)(U) \circ \Pi_1,
\]
where $\Pi_i$ is the regular representation of $(A \otimes \c K) 
\times_{\epsilon_i} G/N$ on $\c H \otimes \ell^2(G) \otimes 
\ell^2(G/N)$ for $i=1,2$ (and $A$ is faithfully represented on a 
Hilbert space $\c H$).  

Finally, from
\cite[Equation (5.1) and Proposition 5.3]{eq:full} 
we have the isomorphism
$\Phi \: A \times_\delta G \times_{\hat\delta} G 
\to A \otimes \c K$ given on generators by
\[
\Phi(a_s,t,r) = a_s \otimes \lambda_s M_{\chi_t} \rho_r,
\]
where $\lambda$ and $\rho$ are the left and right regular
representations of $G$ and $\chi_t$ denotes the characteristic
function of the singleton $\{t\}$.  ($\Phi$ is
the isomorphism of Katayama's duality theorem \cite[Theorem 8]{kat},
but for maximal coactions rather than reduced ones.)
The arguments of \cite{kat}, adapted to our 
context, show that $\Phi$
is equivariant for the coactions $\deltahathat\bigm|$ 
and $\epsilon_1$ of $G/N$; we define $\Theta_1 = \Phi\times G/N$ 
to be the corresponding 
isomorphism of the crossed products.

Careful study of the isomorphisms $\Theta_1$, $\Theta_2$, and 
$\Theta_3$ now shows that the composition 
$\Theta \: A \times_\delta G \times_{\hat\delta}G 
\times_{\deltahathat\bigm|} G/N 
\iso (A \times_{\delta|} G/N) \otimes \c K$ is given on the 
generators by
\[
\Theta(a_s,t,r,gN) = (a_s,trqN) \otimes \lambda_s M_{\chi_t} \rho_r.
\]
Using this, straightforward calculations show that there is an
$ A \times_\delta G \times_{\hat\delta}G
\times_{\deltahathat\bigm|} G/N - A\times_{\delta|}G/N$
imprimitivity bimodule isomorphism
\[
\Upsilon \: Y_{G/G}^G(A) \times_{\delta_Y} G/N
\iso (A \times_{\delta|}G/N) \otimes \ell^2(G)
\]
defined on the generators by
\[
\Upsilon(a_s,t,rN) = (a_s,rN) \otimes \chi_{st}.
\]


\end{proof}

\begin{rem}
Taking $N=G$ in \thmref{Ng-thm} shows that Katayama's bimodule (by
which we mean the bottom and right-hand sides of that rectangle, taken
together) is isomorphic to Mansfield's in this special case.  This
justifies the idea that Mansfield's theorem ``reduces to Katayama's''
when $N=G$, a fact which is well-known to the cognoscenti, but to our
knowledge has not explicitly appeared in the literature.
\end{rem}

To complete the connection with Ng's theorem, we need to pass to 
\emph{reduced} coactions and \emph{amenable} subgroups
in Diagram~\ref{ng diag}.  
Here we use $Y(A,\delta)$ to denote 
the $A \times_\delta G \times_{\hat\delta|}
N$ -- $A \times_{\delta|} G/N$ imprimitivity bimodule 
provided by the original form of Mansfield's Imprimitivity Theorem
\cite[Theorem 27]{man}.
The isomorphism $\Theta_r\colon
A\times_\delta G\times_{\hat\delta}G\times_{\deltahathat|} G/N\to
(A \times_{\delta|} G/N) \otimes \c K(\ell^2)$ is constructed as in
Equation~\eqref{theta-eq}.

\begin{cor}\label{red-Ng-cor}
If $(A,G,\delta)$ is a reduced coaction of a discrete group $G$ and $N$
is an amenable normal subgroup of $G$, then Ng's bimodule is isomorphic
to Mansfield's; that is, the diagram
\[
\xymatrix
{
{A \times_{\delta} G \times_{\hat\delta|} N}
\ar[r]^-{Y(A,\delta)}
&{A \times_{\delta|} G/N}
\\
{A \times_\delta G \times_{\hat\delta} G
\times_{\deltahathat\bigm|} G/N}
\ar[u]^{X(A \times_\delta G)}
\ar[r]^-{\cong}_-{\Theta_r}
&{(A \times_{\delta|} G/N) \otimes \c K(\ell^2)}
\ar[u]_{(A \times_{\delta|} G/N) \otimes \ell^2}
}
\]
commutes.
\end{cor}

\begin{proof}
Since $G$ is discrete, the coaction $\delta$ is automatically
nondegenerate, so by \cite[Theorem 4.7]{qui:fullreduced} there is a unique
full coaction $\delta^f$ of $G$ on $A$ whose reduction coincides with
$\delta$, and then \cite[Proposition 3.8]{qui:fullreduced} gives an
isomorphism $A \times_\delta G \iso A \times_{\delta^f} G$; it is easy
to see that this isomorphism is equivariant for the dual actions.  Then
\cite[Proposition 5.3]{eq:full} applies, giving a
maximal coaction
$(A^m,G,\delta^m)$
(the ``maximalization'' of $\delta^f$)
and an equivariant surjection $\Psi \: A^m \to A$ whose integrated
form $\Psi \times G \: A^m \times_{\delta^m} G \to A \times_{\delta^f}
G$ is an isomorphism which is equivariant for the dual actions.
Then $\Psi$ is also equivariant for the restricted coactions $\delta^m|$
and $\delta^f|$, hence certainly gives a surjection
\[
\Psi \times G/N \: A^m \times_{\delta^m|} G/N
\to A \times_{\delta^f|} G/N.
\]
Since $\delta^f$ is the normalization of $\delta^m$, \cite[Theorem
3.4]{eq:full} tells us that, if $I = \ker \Psi \times G/N$,
then the ideal of 
$A^m \times_{\delta^m} G \times_{\what{\delta^m}|} N$ induced from $I$
via the Mansfield imprimitivity bimodule $Y(A^m)$ coincides with the kernel
of the regular representation
\[
A^m \times_{\delta^m} G \times_{\what{\delta^m}|} N
\to A^m \times_{\delta^m} G \times_{\what{\delta^m}|,r} N,
\]
so that $Y/(Y \cdot I)$ is canonically an $A^m \times_{\delta^m} G
\times_{\what{\delta^m}|,r} N$ -- $A \times_{\delta^f|} G/N$
imprimitivity bimodule. But $N$ is amenable, so the regular
representation of $A^m \times_{\delta^m} G \times_{\what{\delta^m}|} N$
is faithful. Hence we must have $I=\{0\}$, so $\Psi \times G/N$ is
actually an isomorphism of $A^m \times_{\delta^m|} G/N$ onto $A
\times_{\delta^f|} G/N$.

It is now clear from the constructions that the identity map on the
ordered pairs $\{(a_s,t) : s,t \in G\}$ extends to an isomorphism
\[
Y(A^m) \iso Y(A,\delta^f)
\]
of the Mansfield imprimitivity bimodule
$Y(A^m)$ used in the present paper onto the version of the Mansfield
bimodule associated to the normal coaction $\delta^f$ in
\cite{kq:imprimitivity}.
Thus the diagram
\[
\xymatrix
{
{A^m \times_{\delta^m} G \times_{\what{\delta^m}|} N}
\ar[r]^-{Y(A^m)}
\ar[d]_{\cong}
&{A^m \times_{\delta^m|} G/N}
\ar[d]^{\cong}
\\
{A \times_{\delta^f} G \times_{\what{\delta^f}|} N}
\ar[r]_-{Y(A,\delta^f)}
&{A \times_{\delta^f|} G/N}
}
\]
commutes.
On the other hand, one of the main points of \cite{kq:imprimitivity}
is that the diagram
\[
\xymatrix
{
{A \times_{\delta^f} G \times_{\what{\delta^f}|} N}
\ar[r]^-{Y(A,\delta^f)}
\ar[d]_{\cong}
&{A \times_{\delta^f|} G/N}
\ar[d]^{\cong}
\\
{A \times_{\delta} G \times_{\what{\delta}|} N}
\ar[r]_-{Y(A,\delta)}
&{A \times_{\delta|} G/N}
}
\]
commutes; combining these shows that the
top quadrilateral of the diagram
\[
\xymatrix
@C-20pt
{
{A^m \times_{\delta^m} G \times_{\what{\delta^m}|} N}
\ar[rrr]^-{Y(A^m)}
\ar[dr]_{\cong}
&&&{A^m \times_{\delta^m|} G/N}
\ar[dl]^{\cong}
\\
&{A \times_\delta G \times_{\hat\delta|} N}
\ar[r]^-{Y(A,\delta)}
&{A \times_{\delta|} G/N}
\\
&{A \times_\delta G \times_{\hat\delta} G
\times_{\deltahathat\bigm|} G/N}
\ar[u]|{X(A \times_\delta G)}
\ar[r]^-{\cong}_-{\Theta_r}
&{(A \times_{\delta|} G/N) \otimes \c K}
\ar[u]|{(A \times_{\delta|} G/N) \otimes \ell^2}
\\
{A^m \times_{\delta^m} G \times_{\what{\delta^m}}G
\times_{\deltahathat\bigm|} G/N}
\ar[uuu]|{X(A^m \times_{\delta^m} G)}
\ar[ur]_{\cong}
\ar[rrr]^-{\cong}_-{\Theta}
&&&{(A^m \times_{\delta^m|} G/N) \otimes \c K}
\ar[ul]_{\cong}
\ar[uuu]|{(A^m \times_{\delta^m|} G/N) \otimes \ell^2}
}
\]
commutes.
Since the outer rectangle commutes by \thmref{Ng-thm}, and the left,
right, and bottom quadrilaterals are easily seen to commute, we
conclude that the inner rectangle commutes as well. 
\end{proof}


\providecommand{\bysame}{\leavevmode\hbox to3em{\hrulefill}\thinspace}

\end{document}